\documentclass[12pt]{amsart}
\usepackage{amsfonts}

\newtheorem{theorem}{Theorem}

\newtheorem{lemma}[theorem]{Lemma}
\newtheorem{proposition}[theorem]{Proposition}

\def\NN{\hbox{\sf I\kern-.13em\hbox{N}}}
\def\RR{\hbox{\sf I\kern-.14em\hbox{R}}}

\def\Cc{\hbox{\sf C\kern -.47em {\raise .48ex \hbox{$\scriptscriptstyle |$}}
   \kern-.5em {\raise .48ex \hbox{$\scriptscriptstyle |$}} }}

\newcommand{\be}{\begin{equation}}
\newcommand{\ee}{\end{equation}}
\newcommand{\Lp}{L^{\prime}}
\newcommand{\cP}{{\mathcal P}}
\newcommand{\cM}{{\mathcal M}}

\begin{document}

\title{A generalization of Levinger's theorem to positive kernel operators}
\author{Roman Drnov\v{s}ek}
\date \today

\begin{abstract}
We prove some inequalities for the spectral radius of positive operators on 
Banach function spaces. 
In particular, we show the following extension of Levinger's theorem. 
Let $K$ be a positive compact kernel operator on $L^2(X, \mu)$ with the 
spectral radius $r(K)$. Then the function $\phi$ defined by $\phi(t) = r(t K + (1-t) K^*)$ 
is non-decreasing on $[0, \frac{1}{2}]$. We also prove that 
$\| A + B^* \| \ge 2 \cdot \sqrt{r(A B)}$ for any 
positive operators $A$ and $B$ on $L^2(X, \mu)$. 
\end{abstract}

\maketitle

\noindent
{\it Math. Subj.  Classification (2000)}: 47B34, 47B65, 47A10, 47A12, 47A63 \\
{\it Key words}: Banach function spaces, kernel operators, spectral radius,
numerical radius, operator inequalities

\section{Introduction}
\vspace{5mm}

In general there is no relation between the spectral radius of a sum of operators 
on a Banach space and the sum of the corresponding spectral radii. 
So, under appropriate assumptions, any inequality between these two numbers 
might be interesting.
In \cite{Dr92} we proved some inequalities for the spectral radius of a sum of 
positive compact kernel operators on a Banach function space. We thus extended 
the corresponding matrix results proved in \cite{EJ89}. In this article we show
their further generalizations by removing several assumptions from the results in \cite{Dr92}.
As an application of our main result we obtain an extension of Levinger's 
theorem to positive kernel operators on $L^2$-spaces. This beautiful result, 
stated without proof in \cite{Le70}, asserts that 
for a non-negative (square) matrix $A$ the function 
$$ \phi(t) = r(t A + (1-t) A^T) $$
is non-decreasing on $[0, \frac{1}{2}]$ and is non-increasing on 
$[\frac{1}{2}, 1]$. In particular, for all $t \in [0, 1]$, 
the following inequality holds 
$$ r(t A + (1-t) A^T) \ge r(A)  . $$
This theorem was generalized in Bapat \cite{Ba87}, where an elementary proof is given.
Recently, Alpin and Kolotilina \cite[Theorem 7]{AK98} further extended Bapat's result.
Our Theorem \ref{Levinger} includes their extension as a special case.
Finally, Theorem \ref{new} proves an inequality that seems to be new even 
in the finite-dimensional case. For the theory of Banach function spaces and Banach lattices
we refer the reader to the books \cite{Za83}, \cite{Me91} and  \cite{AB85}. 
Here we shall recall some relevant facts.

Let $\mu$ be a $\sigma$-finite positive measure on a $\sigma$-algebra
$\cM$ of subsets of a non-void set $X$.
Let $M(X,\mu)$ be the vector space of all equivalence classes of (almost everywhere equal)
complex measurable functions on $X$. A Banach space $L \subseteq M(X,\mu)$ is
called a {\it Banach function space} if $f \in L$, $g \in M(X,\mu)$,
and $|g| \le |f|$ imply that $g \in L$ and $\|g\| \le \|f\|$.
Throughout the paper, it is assumed that the dimension of $L$ is greater than one and 
that $X$ is the carrier of $L$, that is, there is no subset $Y$ of $X$ of 
strictly positive measure with the property that $f = 0$ a.e. on $Y$ for all $f \in L$
(see \cite{Za83}). The cone of positive elements in $L$ is denoted by $L_{+}$.
A non-negative function $f \in L_{+}$ is said to be {\it strictly positive}
if $f(x) > 0$ for almost all $x \in X$. 
The norm of $L$ is said to be a {\it weakly Fatou norm}  if there exists
a finite constant $k \ge 1$ such that $0 \le f_{\tau} \uparrow f$ in $L$
implies that $\|f\| \le k \cdot \sup_{\tau} \|f_{\tau}\|$.

By $L^{\prime}$ we denote the {\it associate space} (also called the K\"{o}the
dual) of all $g \in M(X,\mu)$ such that
$$ \varphi_g (f) = \int_X \! f \, g \, d\mu $$
defines a bounded linear functional $\varphi_g$ on $L$.
The space $L^{\prime}$ is also a Banach function space with
respect to the associate norm $\| \cdot \|^{\prime}$ defined by
$$ \| g \|^{\prime} = \| \varphi_g \| = 
                      \sup \left\{ \int_X \! |f \, g| \, d\mu : \,
                      f \in L , \, \| f \| \le 1 \right\} ,$$
and it may be considered as a closed subspace of the dual Banach lattice $L^*$.
In view of the definition of $\| \cdot \|^{\prime}$
the following generalized H\"{o}lder's inequality holds
$$ \int_X \! |f \, g| \, d\mu \le \| f \| \, \| g \|^{\prime} $$
for $f \in L$ and $g \in L^{\prime}$.
Note that the set $X$ is also the carrier of the associate space $L^{\prime}$,
and $L^{\prime}$ separates points of $L$ (see \cite[Theorem 112.1]{Za83}).
For any non-negative functions $f$ and $g$ on $X$ we introduce the following notation 
$$ \langle f, g \rangle =  \int_X \! f \, g \, d\mu  . $$
For brevity, the integration over the whole set $X$ will be denoted by 
$\int \, d\mu(x)$ or even $\int \, dx$.

By an {\it operator} on a Banach function space $L$ we always mean a linear
operator on $L$. The spectrum and the spectral radius of a bounded operator $T$ on $L$ are
denoted by $\sigma(T)$ and $r(T)$, respectively. 
An operator $T$ on $L$ is said to be {\it positive} 
if $T f \in L_+$ for all $f \in L_+$. Given operators $S$ and $T$ on $L$,
we write $S \ge T$ if the operator $S - T$ is positive.
It should be recalled that a positive operator $T$ on $L$ is automatically 
bounded and that $r(T)$ belongs to the spectrum of $T$.  
An operator $K$ on $L$ is called a {\it kernel operator} if
there exists a $\mu \times \mu$-measurable function
$k(x,y)$ on $X \times X$ such that, for all $f \in L$ and for almost all $x \in X$,
$$ \int_X |k(x,y) f(y)| \, d\mu(y) < \infty \ \ \ {\rm and} \ \ 
   (K f)(x) = \int_X k(x,y) f(y) \, d\mu(y)  .$$
One can check that a kernel operator $K$ is positive iff 
its kernel $k$ is non-negative almost everywhere. 
We say that $K$ is {\it reducible} if there exists a set $A \in \cM$  such that 
$\mu(A) > 0$, $\mu(A^c) > 0$ and $k=0$ a.e. on $A \times A^c$. Otherwise, if there is
no such set, $K$ is said to be {\it irreducible}.

Let $K$ be a positive kernel operator on $L$ with kernel $k$. 
It is easily seen that $L^{\prime}$ is invariant under the adjoint operator $K^*$.
We denote by $K^{\prime}$ the restriction of $K^*$ to $L^{\prime}$.
One can show \cite[Section 97]{Za83} that $K^{\prime}$ is also a positive kernel operator 
with the kernel $k^{\prime}(x,y) = k(y,x)$ ($x, y \in X$).
The following important observation was already stated in \cite{Dr00} 
for general Banach lattices.

\begin{proposition}
Let $L$ be a Banach function space with a weakly Fatou norm.
If $K$ is a kernel operator on $L$, then $r(K^{\prime}) = r(K)$.
\label{radius_prime}
\end{proposition}

\begin{proof}
It follows from  \cite[Theorem 107.7]{Za83} (see also the equality (2) on
p. 393 of \cite{Za83}) that the space $L$ can be (not necessarily isometrically) embedded 
into $(L^{\prime})^{\prime}$ as a Banach space. Then we have
$r(K) \ge r(K^{\prime}) \ge r((K^{\prime})^{\prime}) \ge r(K)$,
and so $r(K^{\prime}) = r(K)$.
\end{proof}

The following important result is contained in \cite[Theorems 4.13 and 3.14]{Gr95}. 

\begin{theorem} 
Let $K$ be an irreducible positive kernel operator on a Banach function space $L$ 
such that $r(K)$ is a pole of the resolvent $(\lambda - K)^{-1}$.
Then $r(K) > 0$, $r(K)$ is an eigenvalue of $K$ of algebraic multiplicity one, and the 
corresponding eigenspace is spanned by a strictly positive function. 
\label{Jentzsch}
\end{theorem}

It is well-known that the assumption that $r(K)$ is a pole of the resolvent 
$(\lambda - K)^{-1}$ is satisfied if some power of $K$ is a compact operator.
In this case Theorem \ref{Jentzsch} is known as the theorem of Jentzsch and 
Perron (see \cite[Theorem 5.2]{Gr95}).

We will also need the following simple result. 

\begin{proposition} 
Assume that a positive operator $T$ on a Banach function space $L$ is 
the norm limit of a sequence $\{T_n\}_{n \in \NN}$ of positive operators on $L$ such that 
$T_1 \ge T_2 \ge \ldots \ge T$. Then 
$$ r(T) = \lim_{n \rightarrow \infty} r(T_n) . $$
\label{approximation}
\end{proposition} 

\begin{proof}
The sequence $\{r(T_n)\}_{n \in \NN}$ is non-increasing and bounded below by 
$r(T)$, so that $r(T) \le  \lim_{n \rightarrow \infty} r(T_n)$. 
Since the spectral radius is upper semicontinuous, the equality holds in this inequality. 
\end{proof}

\vspace{5mm}

\section{General Banach function spaces}
\vspace{5mm}

Throughout this section, let $L$ be a Banach function space with a weakly Fatou norm.
For brevity, we denote by $L_{++}^{\infty}(X, \mu)$ 
the set of all strictly positive functions $f \in L^{\infty}(X, \mu)_{+}$ 
satisfying $1/f \in L^{\infty}(X, \mu)_{+}$. 
For $d \in L^{\infty}(X, \mu)_{+}$ the {\it multiplication operator} $D$ is a positive 
operator on $L$ defined by $D f = d \, f$. 
Clearly, $D$ is invertible iff $d \in L_{++}^{\infty}(X, \mu)$.

The following lemma that extends \cite[Lemma 2.2]{Dr92} is needed in 
the proof of Theorem \ref{Sum}. 

\begin{lemma}
Let $K$ be a positive kernel operator on $L$ with $r(K) = 1$. Let 
$d$ and $e$ be strictly positive functions in $L_{++}^{\infty}(X, \mu)$, and let 
$D$ and $E$ be the corresponding multiplication operators on $L$.
Let $f \in L_+$ and $g \in \Lp_+$ be strictly positive 
functions such that $K f$ is a strictly positive function satisfying 
$$ {K f \over f} = {K^{\prime} g \over g} \ \ \ {\it and} \ \ \
    \langle K f, g \rangle = 1  . $$
Then 
\be
\langle D K E u, v \rangle \ge 
\exp \left( \int_X \! K \! f \, g  \, \log (d \, e) \, d\mu \right) 
\label{one}
\ee
for any $u \in L_+$ and for any nonnegative measurable function $v$ on $X$ satisfying 
$u \, v = f \, g$.  
If, in addition,  $\langle K u, v \rangle < \infty$, then 
\be
\langle K u, v \rangle \ge  
\exp \left( \int_X \! K \! f \, g \, \log \left({K \! u \over u} \,  {f \over K \! f} \right) 
          \, d\mu \right) \ge 1  . 
\label{two}
\ee
\label{BasicLemma}        
\end{lemma}

\begin{proof}
Since $\langle K f, g \rangle = 1$, the integral in (\ref{one}) exists, 
while it will be seen below that the integral in (\ref{two}) exists provided
$\langle K u, v \rangle < \infty$. 
In fact, there is no loss of generality in assuming that 
$\langle D K E u, v \rangle < \infty$, and consequently, 
$\langle K u, v \rangle < \infty$, since it holds
$$ \langle K u, v \rangle \le \| 1/d \|_{\infty} \cdot \| 1/e \|_{\infty} \cdot 
   \langle D K E u, v \rangle . $$
We will first show the right-hand inequality in (\ref{two}), that is 
\be
\int_X \! K \! f \, g \, \log \left({K \! u \over u} \,  {f \over K \! f} \right) 
          \, d\mu \ge 0 . 
\label{claim}
\ee

We consider the special case when $v \in \Lp_+$. 
For almost all $x \in X$ we define the probability measure on $\cM$ by
$$ \nu_x(A) = {1 \over (Kf)(x)} \int_A k(x,y) \, f(y) \, dy , $$
where $k$ is the kernel of $K$. Using the estimate $|\log(t)| \le t + {1 \over t}$ 
($t > 0$) we obtain that 
\be 
  \int \! K \! f \, g \, \left| \log \left( {u \over f} \right) \right| \, d\mu  \le 
  \int \! K \! f \, g \, \left( {u \over f} + {f \over u} \right) \, d\mu = 
  \langle u , K^\prime g \rangle + \langle K f , v \rangle < \infty . 
\label{L1}  
\ee  
Now, we have 
$$ \int \! (K \! f)(y) \, g(y) \, \log \left( {u(y) \over f(y)} \right) \, dy  = 
\int \! f(y) \, (K^{\prime} \! g)(y) \,  \log \left( {u(y) \over f(y)} \right) \, dy  = $$
$$ =  \int \! f(y) \, \log \left( {u(y) \over f(y)} \right) 
          \left( \int k(x,y) \, g(x) \, dx \right) dy . $$ 
Because of (\ref{L1}) we can use Fubini's theorem to get 
$$ \int \! (K \! f)(y) \, g(y) \, \log \left( {u(y) \over f(y)} \right) \, dy  =  \int \! 
   g(x) \left( \int k(x,y) \, f(y) \, \log \left( {u(y) \over f(y)} \right) dy \right) dx = $$
$$ =  \int \! (K \! f)(x) \, g(x) \,  
   \left( \int \log \left( {u(y) \over f(y)} \right) d\nu_x(y) \right) dx . $$ 
Then, an application of Jensen's inequality gives the inequality
$$ \int \! (K \! f)(y) \, g(y) \, \log \left( {u(y) \over f(y)} \right) \, dy \le 
   \int \! (K \! f)(x) \, g(x) \, 
   \log  \left( \int {u(y) \over f(y)} \, d\nu_x(y) \right) dx = $$ 
$$ = \int \! (K \! f)(x) \, g(x) \, 
   \log  \left( {(K \! u)(x) \over (K \! f)(x)} \right) dx , $$ 
from which (\ref{claim}) follows. 
To prove the general case, define sequences 
$\{u_n\}_{n \in \NN}$ and $\{v_n\}_{n \in \NN}$ 
of strictly positive functions by $u_n = u + f/n$ and $v_n = f g/u_n$. 
Since $v_n \le n \, g$, we have $v_n \in \Lp_+$, and so 
\be
  \int_X \! K \! f \, g \,  
  \log \left({K \! u_n \over u_n} \, {f \over K \! f} \right) \, d\mu \ge 0 ,
\label{three}  
\ee
by the special case of (\ref{claim}). Since 
$$ {K u_n \over u_n} - {K f \over f} = 
   {u \over u_n} \left( {K u \over u} - {K f \over f} \right) , $$
it holds that 
$$ \left\{ x \in X: {(K u_n) (x) \over u_n (x)} \ge {(K f)(x) \over f(x)} \right\} = 
   \left\{ x \in X: {(K u)(x) \over u(x)} \ge {(K f)(x) \over f(x)} \right\} , $$
and the sequence $\left\{ {K u_n \over u_n} \right\}_{n \in \NN}$ is non-decreasing 
on this set. Then, by the Monotone Convergence Theorem, 
\be 
   \lim_{n \rightarrow \infty} 
   \int_X \! K \! f \, g \, \log^+ \left({K u_n \over u_n} \, {f \over K f} \right) \, d\mu = 
   \int_X \! K \! f \, g \, \log^+ \left({K u \over u} \, {f \over K f} \right) \, d\mu  , 
\label{plus}
\ee
where the limit is finite. Namely, using the inequality 
$\log^+ t \le t$ ($t > 0$) we obtain that 
$$ 0 \le \int_X \! Kf \, g  \log^+ \left({K u \over u} {f \over K f} \right) \, d\mu \le 
   \int_X \! Kf \, g  \, {K u \over u} \, {f \over K f} \, d\mu = 
   \langle K u , v \rangle < \infty  . $$
This shows that the integral in (\ref{claim}) is defined
(and its value belongs to $[- \infty, \infty)$). Similarly, we obtain that 
$$ \lim_{n \rightarrow \infty} 
   \int_X \! K \! f \, g \, \log^- \left({K u_n \over u_n} \, {f \over K f} \right) \, d\mu = 
   \int_X \! K \! f \, g \, \log^- \left({K u \over u} \, {f \over K f} \right) \, d\mu  ,  $$
which together with (\ref{plus}) gives that 
$$ \lim_{n \rightarrow \infty} 
   \int_X \! K \! f \, g \, \log \left({K u_n \over u_n} \, {f \over K f} \right) \, d\mu = 
   \int_X \! K \! f \, g \, \log \left({K u \over u} \, {f \over K f} \right) \, d\mu .  $$
In view of (\ref{three}) this completes the proof of (\ref{claim}). 

We now define the probability measure $\lambda$ on $\cM$ by
$$ \lambda (A) = \int_A K \! f \, g \, d\mu . $$
An application of Jensen's inequality gives that 
$$ \log \left( \langle K u, v \rangle \right) = 
   \log \left( \int {Ku \over u} \,  {f \over Kf}  \, d\lambda \right) \ge  $$
$$ \ge \int \log \left( {Ku \over u} \, {f \over Kf}  \right) d\lambda = 
   \int \! K \! f \, g \, \log \left({K u \over u} \, {f \over K f} \right) \, d\mu , $$
so that the left-hand inequality holds in (\ref{two}). 
Similarly, we have 
$$ \log \left( \langle D K E u, v \rangle \right) = 
   \log \left( \int d \, e \, {K (E u) \over E u} \, {f \over Kf}  \, d\lambda \right) \ge  $$
$$ \ge \int \log \left( d \, e \, {K (E u) \over E u} \, {f \over Kf}  \right) d\lambda = 
   \int \! K \! f \, g  \, \log (d \, e) \, d\mu +  \int \! K \! f \, g \, \log 
   \left({K (E u) \over E u} \, {f \over K f} \right) \, d\mu . $$
Since the last integral is non-negative by (\ref{claim}), this gives (\ref{one}). 
\end{proof}

The following result extends Theorems 2.4 and 2.6 in \cite{Dr92}. Its 
finite-dimensional version was shown in \cite[Theorem 2.3]{EJ89}.

\begin{theorem}
Let $K_1$, $K_2$, $\ldots$, $K_n$ be positive kernel operators on $L$. 
Assume that $f_1$, $f_2$, $\ldots$, $f_n \in L_+$ and $g_1$, $g_2$, $\ldots$, 
$g_n \in \Lp_+$ are strictly positive functions satisfying 
$$ K_i f_i = r(K_i) f_i \ , \ \ K_i^{\prime} g_i = r(K_i) g_i $$
and be normalized so that 
$$ f_i \cdot g_i = h \ (i = 1,2, \ldots, n) \ \ \ {\rm and} \ \ \ 
   \int_X \! h \, d\mu = 1 .$$
Furthermore, let $d_1$, $\ldots$, $d_n$ and $e_1$, $\ldots$, $e_n$ be in 
$L^{\infty}(X, \mu)_+$,  and 
let $D_1$, $\ldots$, $D_n$ and  $E_1$, $\ldots$, $E_n$ be the 
corresponding multiplication operators on $L$. Then 
\be 
r(\sum_{i=1}^n D_i K_i E_i) \ge \sum_{i=1}^n r(K_i) \, 
\exp \left( \int_X \! h \log (d_i e_i) \, d\mu \right) 
\label{main_inequality}
\ee
adopting the convention $\exp (-\infty) = 0$. 
In particular, for all positive numbers $t_1$,  $\ldots$, $t_n$, 
\be 
r(t_1 K_1 + \ldots + t_n K_n) \ge t_1 r(K_1) + \ldots + t_n r(K_n) .
\label{sum_inequality}
\ee
\label{Sum}
\end{theorem}

\begin{proof}
If, for some $i$, $d_i e_i = 0$ on the set of positive measure, then 
$\int_X \! h \log (d_i e_i) \, d\mu = - \infty$, which together with the monotonicity 
of the spectral radius convinces us that 
there is no loss of generality in assuming that $\{d_i\}_{i=1}^n$ and $\{e_i\}_{i=1}^n$
are strictly positive functions. Also, we may assume that $r(K_i) > 0$ for all $i$.

Consider first the case when $\{d_i\}_{i=1}^n$ and $\{e_i\}_{i=1}^n$ are in 
$L_{++}^{\infty}(X, \mu)$. 
Denote $K = D_1 K_1 E_1 + \ldots + D_n K_n E_n$, pick $\lambda > r(K)$, and set  
$$ u = (\lambda - K)^{-1} f_1 = \sum_{j=0}^{\infty} {\lambda}^{-j-1} K^j f_1 . $$ 
Then $u$ is a strictly positive function in $L$ satisfying $K u \le \lambda u$. 
Denoting $v = h/u$ we apply (\ref{one}) of Lemma \ref{BasicLemma} 
for the operator $K_i/r(K_i)$, $i=1, \ldots, n$, to get 
$$ \langle D_i K_i E_i u, v \rangle \ge 
   r(K_i) \, \exp \left( \int_X \! h \, \log (d_i e_i) \, d\mu \right) . $$
Summing over $i$ gives the inequality 
$$ \sum_{i=1}^n r(K_i) \, \exp \left( \int_X \! h \, \log (d_i e_i) \, d\mu \right) \le $$
$$ \le \sum_{i=1}^n \langle D_i K_i E_i u, v \rangle = 
   \langle K u, v \rangle \le \langle \lambda u, v \rangle = \lambda  . $$
Since this is true for any $\lambda > r(K)$, the inequality (\ref{main_inequality}) follows.

To remove the assumptions on $\{d_i\}_{i=1}^n$ and $\{e_i\}_{i=1}^n$, 
define $d_i^{(m)} = \max \{d_i, {1 \over m}\}$ and $e_i^{(m)} = \max \{e_i, {1 \over m}\}$
($m \in \NN$, $i=1, \ldots, n$), and let $D_i^{(m)}$ 
and $E_i^{(m)}$ be the corresponding multiplication 
operators on $L$. Then, by the above,  
$$ r(\sum_{i=1}^n D_i^{(m)} K_i E_i^{(m)}) \ge \sum_{i=1}^n r(K_i) \, 
\exp \left( \int_X \! h \log (d_i^{(m)} e_i^{(m)}) \, d\mu \right) . $$
When $m$ tends to infinity, the left-hand side approaches $r(K)$ by Proposition
\ref{approximation}, while 
$$ \lim_{m \rightarrow \infty} \int_X \! h \log (d_i^{(m)} e_i^{(m)}) \, d\mu =
   \int_X \! h \log (d_i e_i) \, d\mu $$
by the Monotone Convergence Theorem (for decreasing sequences). This yields 
the inequality (\ref{main_inequality}), and the proof is finished.
\end{proof}

A glance at the proof above shows that Theorem \ref{Sum} also holds in the case when 
some operators of $K_1$, $K_2$, $\ldots$, $K_n$ are positive multiples of the 
identity operator, or in other words, every $K_i$ is a sum of a positive kernel operator
and a non-negative multiple of the identity.

Given a positive operator $T$ on $L$, let $\cP_+(T)$ denote the set of all 
functions $p(z)=\sum_{k=0}^{\infty} a_k z^k$ such that $a_k \geq 0$ for all $k$ and 
the convergence radius of $p$ is greater than $r(T)$. 
Using the spectral mapping theorem one can show easily
that $r(p(T)) = p(r(T))$ for all $p \in \cP_+(T)$.

\begin{theorem}
Under the assumptions of Theorem \ref{Sum}, let $p_i \in \cP_+(T_i)$ for 
$i=1, \ldots, n$. Then 
$$ r(p_1 (K_1) + \ldots + p_n (K_n)) \ge p_1 (r(K_1)) + \ldots + p_n (r(K_n)) . $$
In particular, if $s_i > r(K_i)$ for $i=1, \ldots, n$, then 
$$ r((s_1 - K_1)^{-1} + \ldots + (s_n - K_n)^{-1}) \ge 
   {1 \over s_1 - r(K_1)} + \ldots + {1 \over s_n - r(K_n)} . $$
\end{theorem}

\begin{proof}
We first claim that every $p_i(K_i)$, $i=1, \ldots, n$, is the sum of a kernel operator 
and a non-negative multiple of the identity operator $I$. If 
$p_i(z)=\sum_{k=0}^{\infty} a_k z^k$ with $a_k \geq 0$, then 
$p_i(K_i) - a_0 I$ is the limit (in norm and in order) of an increasing sequence 
of kernel operators. It follows that it is a kernel operator 
(see e.g. \cite[Theorem 94.5]{Za83}). This proves our claim. 
Now, according to the remark following the proof of Theorem \ref{Sum} 
we may apply the inequality (\ref{sum_inequality}) of Theorem \ref{Sum} 
for operators $p_1(K_1)$, $\ldots$, $p_n(K_n)$ to get
$$ r(p_1 (K_1) + \ldots + p_n (K_n)) \ge r(p_1(K_1)) + \ldots + r(p_n(K_n)) = 
p_1 (r(K_1)) + \ldots + p_n (r(K_n)) . $$
\end{proof} 

As an extension of Theorem 4.2 in \cite{FK75} we now show that the inequality  
(\ref{main_inequality}) of Theorem \ref{Sum} for $n=1$ can be improved if the 
operator is of the form $(s - K)^{-1}$, where $s > r(K)$.

\begin{theorem}
Let $K$ be a positive operator on $L$ with $r(K) > 0$ that is a sum of a 
positive kernel operator and a non-negative multiple of the identity. 
Assume that $f \in L_+$ and $g \in \Lp_+$ are strictly positive functions satisfying 
$K f = r(K) f$, $K^{\prime} g = r(K) g$ and $\langle f, g \rangle = 1$. 
Let $d$ be in $L^{\infty}(X, \mu)_+$, and let $D$ be the corresponding 
multiplication operator on $L$. Then 
\be 
r(D K) \ge r(K) \, \exp \left( \int_X \! f \, g \, \log (d) \, d\mu \right) . 
\label{weaker}
\ee
Furthermore, for $s > r(K)$ it holds 
\be 
r(D (s - K)^{-1}) \ge r((s - K)^{-1}) \, \left( \int_X \! f \, g \, d \, d\mu \right) . 
\label{stronger} 
\ee
\end{theorem}

\begin{proof}
The inequality (\ref{weaker}) is a special case of (\ref{main_inequality}).
Denote $T = (s - K)^{-1}$ and pick $\lambda > r(D T)$. 
Then $w = (\lambda - D T)^{-1} f$ is a strictly positive function in $L$ 
satisfying $D T w \le \lambda w$. 
Set $u = T w$ and $v = f \cdot g / u$. If we apply (\ref{two}) of Lemma \ref{BasicLemma}        
for the operator $K/r(K)$, we obtain that 
$\langle K u, v \rangle \ge r(K)$, and so 
$$ \langle T^{-1} u, v \rangle = \langle (s - K)u, v \rangle \le 
    s - r(K) = {1 \over r(T)} . $$
On the other hand, since $\lambda T^{-1} u = \lambda w \ge D T w = d u$, we have 
$\lambda \langle T^{-1} u, v \rangle \ge \langle d u, v \rangle$. It follows that 
$\lambda \ge r(T)  \langle d u, v \rangle$ which implies (\ref{stronger}). 
\end{proof}

Observe that (\ref{stronger}) is really a sharpening of  (\ref{weaker}) 
for the special class of positive operators, since 
$$ \exp \left( \int_X \! f \, g \, \log (d) \, d\mu \right) \le 
   \int_X \! f \, g \, d \, d\mu $$
by Jensen's inequality. Also, simple examples show that in (\ref{weaker}) 
$\exp \left( \int_X \! f \, g \, \log (d) \, d\mu \right)$ can not be replaced 
by $\int_X \! f \, g \, d \, d\mu$. 
(Consider $K =  \left[ \begin{array}{cc}
                      0    &    1  \\
                      1    &    0  
                   \end{array} \right]$ on $L = \Cc^2$.) \\

\section{$L^2$-spaces}
\vspace{5mm}

In \cite{Dr92} we proved an extension of Levinger's 
inequality to positive kernel operators on $L^2$-spaces.
Unfortunately, we were able to show it only under some assumptions 
on the kernel of the operator. We now show that these assumptions are redundant,
as we expected. In the finite-dimensional case this result was proved in 
\cite[Theorem 7]{AK98}. 

\begin{theorem}
Let $K$ be a positive kernel operator on $L^2(X, \mu)$ such that 
$r(K)$ is an isolated point of $\sigma(K)$ and the corresponding 
Riesz idempotent has finite rank. 
Let $d \in L_{++}^{\infty}(X, \mu)$ be a strictly positive function, 
and let $D$ be the corresponding multiplication operator on $L^2(X, \mu)$. 
Then, for any $t \in [0,1]$, 
\be 
r(t D K D^{-1} + (1-t) K^*) \ge r(K) . 
\label{LevEq}
\ee
If, in addition, the operator $K$ is compact and if $\phi : [0,1] \rightarrow [0, \infty)$ 
is defined by 
$$ \phi(t) = r(t D K D^{-1} + (1-t) K^*) , $$
then $\phi$ is non-decreasing on $[0, {1 \over 2}]$ and is non-increasing on 
$[{1 \over 2}, 1]$.
\label{Levinger}
\end{theorem}

\begin{proof}
Consider first the case when $D = I$, the identity on $L$. 
If $K$ is irreducible, then by Theorem \ref{Jentzsch} there exist 
strictly positive functions  $f, g \in L^2(X, \mu)$ satisfying 
$K f = r(K) f$, $K^* g = r(K) g$ and $\langle f, g \rangle = 1$, 
and  the inequality (\ref{LevEq}) follows from Theorem \ref{Sum} with $K_1 = K$, $K_2 = K^*$, 
$f_1 = g_2 = f$ and  $g_1 = f_2 = g$. 
For general $K$ pick any strictly positive function $u \in L^2(X, \mu)$. 
(Such functions exist because the measure $\mu$ is $\sigma$-finite.)
Denote by $K_0$ an irreducible kernel operator with strictly positive kernel 
$u(x) u(y)$ ($x, y \in X$). For each $m \in \NN$ define an irreducible positive kernel 
operator on $L^2(X, \mu)$ by $K_m = K + {1 \over m} K_0$. 
Then $r(K_m) \ge r(K)$, and the left (and, similarly, the right) essential spectra 
of $K_m$ and $K$ coincide.
Now, Proposition XI.6.9 and Theorem XI.6.8 of \cite{Co90} imply that $r(K_m)$ is an isolated 
point of $\sigma(K_m)$ and the corresponding Riesz idempotent has finite rank. 
By the first part of the proof, we then have 
$$ r(t K + (1-t) K^* + {1 \over m} K_0) = 
   r( t(K + {1 \over m} K_0) + (1-t) (K + {1 \over m} K_0)^*)  
   \ge r(K+{1 \over m} K_0) . $$
Letting $m \rightarrow \infty$ we get $\phi(t) \ge r(K)$ by 
Proposition \ref{approximation}, which proves (\ref{LevEq}) in the case $D = I$. 
Since $\phi(t) = \phi(1-t)$, it remains to show in this special case 
that $\phi$ is non-decreasing on $[0, {1 \over 2}]$ provided $K$ is compact.
Let $0 \le t < s \le {1 \over 2}$. Then, by (\ref{LevEq}), 
$$ \phi(t) \le r(u(t K + (1-t) K^*) + (1-u)(t K + (1-t) K^*)^*) = 
 r((2 u t - u - t + 1)K + (t + u - 2 u t)K^*) $$
for all $u \in [0,1]$. Put $u = {1 - s - t \over 1 - 2 t}$ to obtain that 
$\phi(t) \le r( s K + (1 - s) K^*) = \phi(s)$.

The general case follows from the special one. To show this, 
let $E$ be the multiplication operator on $L$ the multiplier of which is $\sqrt{d}$, so that 
$E^2 = D$. Introducing the notation 
$$ \phi_{K, D}(t) = r(t D K D^{-1} + (1-t) K^*)  $$
we have, for all $t \in [0,1]$,
$$ \phi_{K, D}(t) = r(E( t E K E^{-1} + (1-t) E^{-1} K^* E) E^{-1}) = $$
$$ =  r(t E K E^{-1} + (1-t) (E K E^{-1})^*) = \phi_{E K E^{-1}, I}(t) . $$
Since $\phi_{E K E^{-1}, I}(t) \ge r(E K E^{-1}) = r(K)$ by the special case, 
(\ref{LevEq}) follows.
If, in addition, $K$ is compact, then $\phi_{E K E^{-1}, I}$ is non-decreasing 
on $[0, {1 \over 2}]$ and is non-increasing on $[{1 \over 2}, 1]$ by the special case, 
and so the same is also true for $\phi_{K, D}$. This completes the proof.
\end{proof}

We do not know whether Theorem \ref{Levinger} is valid for every positive operator $K$ on 
$L^2(X, \mu)$. However, we shall show below that for $t = 1/2$ the inequality 
(\ref{LevEq}) holds for all  positive operators on $L^2(X, \mu)$.
To do this, we recall that the {\it numerical radius} $w(A)$ of a bounded operator $A$ on 
$L^2(X, \mu)$ is defined by 
$$ w(A) = \sup \{ | \langle A f, f \rangle | : f \in L^2(X, \mu), \| f \|_2 = 1 \} . $$
If, in addition, $A$ is positive, then we have 
$$ w(A) = \sup \{ \langle A f, f \rangle  : f \in L^2(X, \mu)_+ , \| f \|_2 = 1 \} . $$
Indeed, this follows from the estimate
$$  | \langle A f, f \rangle | \le \int_X \! |A f| \, |f| \, d\mu \le 
    \langle A |f|, |f| \rangle  $$
that holds for any  $f \in L^2(X, \mu)$. It is well-known \cite{GD97} that 
$$ r(A) \le w(A) \le \|A\| $$
for all bounded operators $A$ on $L^2(X, \mu)$.  

\begin{theorem}
Let $A$ be a positive operator on $L^2(X, \mu)$. Then, for any $t \in [0,1]$, 
\be 
\|A\| \ge \|t A + (1-t) A^*\| \ge w(t A + (1-t) A^*) = w(A) \ge r(A)
\label{linear}
\ee
and 
\be
\|(t A + (1-t) A^*)^2\| \ge w((t A + (1-t) A^*)^2) \ge  w(A^2) \ge (r(A))^2 .
\label{square} 
\ee
Furthermore, if $d$ is in  $L_{++}^{\infty}(X, \mu)$ and $D$ is the corresponding 
multiplication operator on $L^2(X, \mu)$, then  
\be 
r(D A D^{-1} + A^*) \ge 2 \, r(A) 
\label{symmetric}
\ee
\end{theorem}

\begin{proof}
The equality in (\ref{linear}) follows from 
$$ \langle (t A + (1-t) A^*) f, f \rangle = t \langle A f, f \rangle + 
(1-t) \langle f, A f \rangle =  \langle A f, f \rangle , $$
which holds for all $f \in L^2(X, \mu)_+$. 
The remaining inequalities in (\ref{linear}) are clear. 
Similarly, only the second inequality in (\ref{square}) needs a proof. 
This relation is a consequence of the following inequality 
$$ \langle (t A + (1-t) A^*)^2 f, f \rangle \ge \langle A^2 f, f \rangle $$
that holds for every $f \in L^2(X, \mu)_+$, since it is equivalent to 
$t (1-t) \| A f - A^* f\|^2_2 \ge 0$.  
Setting $t = 1/2$ in (\ref{linear}) we obtain (\ref{symmetric}) in the case $D = I$, since 
$r(A + A^*) = w(A + A^*) = \|A + A^*\|$. 
The general case can be obtained from the special one as in the proof of Theorem 
\ref{Levinger}. Namely, 
if $E$ is the multiplication operator on $L$ with the multiplier $\sqrt{d}$, then 
$$ r(D A D^{-1} + A^*) = r(E( E A E^{-1} + E^{-1} A^* E) E^{-1}) = $$
$$ =  r(E A E^{-1} + (E A E^{-1})^*) \ge  2 \, r(E A E^{-1}) = 2 \, r(A). $$
\end{proof}

An application of Berberian's trick concerning $2 \times 2$ operator matrices gives the 
following result which seems to be new even in the finite-dimensional case.

\begin{theorem}
Let $A$ and $B$ be positive operators on $L^2(X, \mu)$. Then 
$$ \| A + B^* \| \ge 2 \cdot \sqrt{r(A B)} . $$
If, in addition, $A$ and $B$ are compact kernel operators, then, for each $t \in [0,1]$, 
$$ \max\{ \| t A + (1-t) B^* \|, \| t B + (1-t) A^* \|\}  \ge \sqrt{r(A B)} . $$

\label{new}
\end{theorem}

\begin{proof}
Let $T$ be a positive operator on $L^2(X, \mu) \oplus L^2(X, \mu)$ defined by
$2 \times 2$ operator matrix 
$$ T = \left[ \begin{array}{cc} 
                 0 & A \\
                 B & 0 \\
              \end{array} \right] . $$ 
Then $r(T + T^*) = \| T + T^* \| = \| A + B^* \|$ and 
$(r(T))^2 = r(T^2) = r(A B)$. By (\ref{symmetric}), we obtain that 
$$ \| A + B^* \| = r(T + T^*) \ge  2 \, r(T) = 2 \, \sqrt{r(A B)} . $$
If, in addition, $A$ and $B$ are compact kernel operators, then $T$ 
is a compact kernel operator as well. Then, for each $t \in [0,1]$, 
$$ \sqrt{r(A B)} = r(T) \le r(t T + (1-t) T^*) \le $$
$$ \le \| t T + (1-t) T^*\| = \max\{ \| t A + (1-t) B^* \|, \| t B + (1-t) A^* \|\} , $$
where we have used (\ref{LevEq}). This completes the proof.
\end{proof}

\vspace{3mm}

{\it Acknowledgment.} This work was supported in part 
by the Ministry of Education, Science and Sport of Slovenia. \\

\vspace{3mm}

\noindent
Roman Drnov\v sek \\
Department of Mathematics \\
University of Ljubljana \\
Jadranska 19 \\
SI-1000 Ljubljana \\
Slovenia \\
e-mail : roman.drnovsek@fmf.uni-lj.si 

\end{document}